\documentclass[11pt]{article}
\usepackage{amsmath,amssymb,amsfonts,amsthm,graphics,graphicx,psfrag}
\newtheorem{example}{Example}[section]

\newtheorem{lemma}[example]{Lemma}

\def\sign{\,{\rm sign}\,}
\def\eps{\varepsilon}
\def\vphi{\varphi}

\def \F{\mathcal F}
\def \N{\mathbb N}
\def \R{\mathbb R}

\def \T{\mathbb T}
\def \H{{H^1(\R)}}
\def\dis{\displaystyle}
\def\meas{{\rm \,meas\,}}
\def\st{{\rm \ \ s.t.\ \ }}
\begin{document}
\bibliographystyle{plain}
\title {\bf Conservative solution of the Camassa Holm Equation on the real line}
\author
       {Massimo Fonte
        \\
        \small{S.I.S.S.A., Via Beirut, 2/4}
        \\
        \small{34014 Trieste, Italy}
        \\
        \small{e-mail: \tt {fonte@sissa.it}}
       }
\date{}
\maketitle

\begin{abstract}
In this paper we construct a global, continuous flow of solutions to the Camassa-Holm equation on the space $H^1(\R)$. In a previous paper \cite{BF2}, A. Bressan and the author constructed spatially periodic solutions, whereas in this paper the solutions are defined in all the real line. We introduce a distance functional, defined in terms of an optimal transportation problem, which allows us to study the continuous dependance w.r.t. the inital data with a certain decay at infinity.
\end{abstract}

\section{Introduction}
In \cite{CH} the authors present a nonlinear partial differential equation which describes the behaviour of the shallow water as a completely integrable hamiltonian system
\begin{equation}
\label{kdiv0}
u_t+2\kappa u_x - u_{xxt}+3u u_x = 2 u_x u_{xx}+ u u_{xxx}.
\end{equation}
where $u$ is the fluid velocity in the $x$ direction and $\kappa$ is a constant related to the critical shallow-water  wave speed. For the physical description of such equation, we refer to \cite{CH, Jo, CM1, CM2} and to the bibliographic references of \cite{BF2}. In the present paper we study the limit case $\kappa= 0$, a condition in which, starting by a initial smooth data it can develop to a peaked solution with one or more cuspids, the so called multi-peakon function. As an example, in \cite{CHH} a simulation shows that starting from a parabolic initial data in a periodic domain, the system evolves in a train of positive peakons. 
The equation can be written in nonlocal form as a scalar conservation law with an integro-differential source term:
\begin{equation}
u_t+\left(\frac{u^2}{2}\right)_x+P^u_x=0
\label{prblCH}
\end{equation}
where $P^u$ is defined in term of a convolution:
$$
P^u\doteq \frac 12 e^{-|x|}*\left(u^2+\frac{u_x^2}2\right)\,.
$$
Observe that the function $\frac 12 e^{-|x|}$ is the distributional solution of the equation 
$$
\left(Id-\partial_{xx}\right)f=\delta_0
$$
where $\delta_0$ is the Dirac measure centered at the origin.
A \emph{multi-peakon} is a function of the form  
$$
u(x)=\sum_{i=1}^N p_i e^{-|x-q_i|}
$$
it is well known (see \cite{CE1,CHL,HR}) that if such a function is subject to (\ref{prblCH}), its evolution remains of the same shape, and as long as they are well defined, the coefficients $p_i$, $q_i$ are the solution to the system of ODE
$$
\left\{
\begin{array}{l}
\dis \dot q_i=\sum_{j=1}^N p_j e^{-|q_i-q_j|}\,,
\\
\dis \dot p_i=p_i\sum_{j=1}^N p_j \sign(q_i-q_j)e^{-|q_i-q_j|},
\end{array}
\right.\qquad i=1,\dots, N.
$$
In the smooth case, differentiating (\ref{prblCH}) w.r.t. $x$ and multiplying by $u_x$ we obtain the conservation law with source term
\begin{equation}
\label{derevolution}
\left(\frac{u_x^2}{2}\right)_t+\left(\frac{u u_x^2}{2}-\frac {u^3}{3}\right)_x=-u_x\,P.
\end{equation}
The previous equation, together with the \emph{Camassa-Holm} equation (\ref{prblCH}), prove that the total energy
$$
E\doteq\int_{\R}[u^2(t,x)+u_x^2(t,x)]dx 
$$
is a conserved quantity as long as the solution remains regular. Constantin and Escher \cite[Theorem 4.1]{CE1} shown that even if the initial data is sufficiently regular, blow-up of the gradient $u_x$ can occur in finite time whenever the initial data has a negative slope.
In Section \ref{ODEsystem} we implement a technique based on appropriate rescaled variables in order to resolve the singularities which occurs 
at the times where the gradient $u_x$ blows up. The new system of ODE can be solved in a unique way in a neighborhood of the time of blow up and the solution turns out to preserve the energy $E$ also after this time. 

Motivated by the existence of the multi-peakon solutions, whose decay at infinitive is like $e^{-|x|}$, we introduce the space $X_\alpha$ of the $H^1$ function with an exponential decay: let  $0<\alpha<1$, we define 
\begin{equation}
X_\alpha= \left\{u\in H^1(\R) \,:\, C^{\alpha,u}\doteq \int_R \left[u^2(x)+u_x^2(x)\right]e^{\alpha|x|}\,dx<\infty\right\}\,.
\label{decayHP}
\end{equation}
In this space we define a distance that is related to an optimal transportation problem (see \cite{V}). Fetching the theory developed by Bressan and Constantin \cite{BC1} for the \emph{Hunter-Saxton} equation, and by Bressan and Fonte \cite{BF2} for the periodic solution of the Camassa-Holm equation, the topology induced by functional $J$ constructed in Section \ref{metricsection} turns out to be weaker than the $H^1$ topology, but useful because with this metric we prove the stability of the multi-peakon solutions w.r.t. the initial conditions, as we will show in Section \ref{stabcor}. 

\section{Multipeakon solutions}\label{ODEsystem}
In this section we shall construct a solution of the Camassa Holm equation starting from an initial condition $\bar u^\varepsilon\in X_\alpha$ of the form
$$
u_0^\varepsilon(x)=\sum_{j=1}^{N_\varepsilon} p_j e^{-|x-q_j|}\,.
$$
The motivation of this choice is given by the form of traveling wave solution (see \cite{CH} and \cite[Example 5.2]{CE1}).  Looking for solution of the equation (\ref{kdiv0}) in the traveling wave form $u(t,x)=U(x-ct)$, with a function $U$ that vanishes at infinity, the limit of $\kappa\to 0$ leads to the function $U= c e^{-|x-ct|}$. The evolution of an initial data like $u_0^\varepsilon$ remains then of the same shape
$$
u^\varepsilon(t,x)=\sum_{j=1}^{N_\varepsilon} p_j(t) e^{-|x-q_j(t)|}\,.
$$
As long as the classical solution of the problem 
\begin{equation}
\label{HSYS}
\left\{
\begin{array}{l}
\dis \dot q_i=\sum_{j=1}^N p_j e^{-|q_i-q_j|}\,,
\\
\dis \dot p_i=p_i\sum_{j=1}^N p_j \sign(q_i-q_j)e^{-|q_i-q_j|}
\end{array}
\right.
\end{equation}
exists, the solution of this system gives the coefficients ${\bf p }(t)=(p_1,\dots, p_{N_\varepsilon})$ and ${\bf q}(t)=(q_1,\dots, q_{N_\varepsilon})$ for the solution $u^\varepsilon(t,x)$ to the Camassa-Holm equation. 
Observe that the previous system can be viewed as an Hamiltonian system with Hamiltonian function $H({\bf q},{\bf p })=\frac 12~\sum_{i,j}p_i p_je^{-|q_i-q_j|}$.
 
In \cite{HR} the autors prove the existence of a global multi-peakon solution when strengths $p_i$ are positive for all $i=1\dots N_\varepsilon$ and the convergence of the sequence of multi-peakon solution. If $u_0$ is an initial data such that the distribution $u_0-{u_0}_{xx}$ is a \emph{positive} Radon measure, there exists a sequence of multi-peakons that converges in $L^\infty(\R,H^1_{loc}(\R))$. In this case the crucial fact is that no interaction between the peakons occurs, and then the gradient remains bounded. However, a general initial data contains both positive and negative peakons, as in the so called peakon-antipeakon interaction: one positive peakon with strength $p$, centered in $-q$, moves forward and one negative anti-peakon in $q$, with strength $p$ moves backward. The evolution of the system produces the overlapping of the two peakons at finite time $t=\tau$, so that $q\to 0$.  The conservation of the energy $E=H({\bf q}(t),{\bf p }(t))$ yields
\begin{equation}
\label{zetalimit}
E=\lim_{t\to \tau^-} p^2(1-e^{-2|q|})\,.
\end{equation}
and then the quantity $p$ diverges in finite time. At the point $(\tau,0)$ occurs thus a singularity for the solution $u$.  To extend the solution also after the interaction time with a solution which conserves the energy $E$ we can think that at the interaction point emerge an antipeakon/peakon couple, the first, negative, moving backward and the second, positive, moving forward with coefficients $(-q, -p)$ and $(q,p)$. According to the conservation of the energy, the choice of $q$ and $p$ must satisfy (\ref{zetalimit}) as $t\to \tau^+$.
It  yields a change of variables which resolves the singularity at $(\tau,0)$
$$
\zeta\doteq p^2 q\qquad \omega\doteq\arctan (p)
$$ 
with this choice, the Hamiltonian system leads to the ODE
$$
\frac d{dt}
\left(
\begin{array}{c}
\zeta
\\
\omega
\end{array}
\right)
=f(\zeta,\omega)\,,
\qquad
\left(
\begin{array}{c}
\zeta
\\
\omega
\end{array}
\right)(\tau)
=
\left(
\begin{array}{c}
\frac E 2
\\
\frac \pi 2
\end{array}
\right)
$$
with 
$$
f(\zeta,\omega)= 
\left(
\begin{array}{c}
\left[1- e^{-\zeta \cot^2(\omega)}-\zeta\cot^2(\omega)e^{-\zeta \cot^2(\omega)}\right]\tan^3(\omega)
\\
\sin^2(\omega)e^{-\zeta \cot^2(\omega)}
\end{array}
\right)
$$
and $f$ is a Lipschitz vector field in a neighborhood of the point $( \frac E 2,\frac \pi 2)$. The solution $(\zeta(t),\omega(t))$ of this problem provides then the unique couple $(q(t),p(t))$ which coincides with the classical solution of the Hamiltonian system for $t<\tau$ and extends it for $t\geq \tau$.

This example suggests the way to construct the multi-peakon solution whenever an interaction between peakons occurs (see also \cite{BF2} for an ``energetic'' motivation). Suppose that two or more peakons with strengths $p_1,\dots,p_k$ annihilate at the position $\bar q$ at time  $\tau$  and produce a blow up of the gradient $u_x$. The conservation of the energy yields that there exists and is positive the limit 
$$
e_{\tau}\doteq \lim_{t\to \tau^-} \int_{\xi^-(t)}^{\xi^+(t)} u_x^2(t,x)\,dx
$$ 
where $\xi^-$ and $\xi^+$ are the smallest and the largest characteristic curve passing through the point $(\tau,\bar q)$. Assume that after the interaction appear two peakons with strengths $p_1,\,p_2$ and placed at the position $q_1,\,q_2$.
Let consider the change of variables
$$
z=p_2+p_1\quad w=2\arctan (p_2-p_1)\quad \eta = q_2+q_1 \quad \zeta=(p_2-p_1)^2(q_2-q_1),
$$
then the system (\ref{HSYS}) turns out to be
$$
\begin{array}{rl}
\dot w=&\!\! -\left[\sin (w) \cosh \Big(\frac {\zeta}{2\tan^2 (w/2)}\Big)
+2z \sinh\Big(\frac {\zeta}{2\tan^2 (w/2)}\Big) \right]\cdot\sum\limits_{j\geq k+1} p_j e^{-q_j+ \eta /2}
\\
&+[z^2 \cos^2(w/2) - \sin^2 (w/2)]e^{-\frac{\zeta}{\tan^2 (w/2)}}
\\
\dot z= &\!\!-\left[\frac 12\sin(w) \sinh \Big(\frac {\zeta}{2\tan^2 (w/2)}\Big)+ z\cosh \Big(\frac {\zeta}{2\tan^2 (w/2)}\Big) \right]\cdot \sum\limits_{j\geq k+1} p_j e^{-q_j}
\\
\dot \eta= &\!\! z[1+ e^{-\frac{\zeta}{\tan^2 (w/2)}}]+ 2 \cosh \Big(\frac {\zeta}{2\tan^2 (w/2)}\Big)\cdot \sum\limits_{j\geq k+1} p_j e^{-q_j}
\\
\dot \zeta = &\!\!\frac{z^2\zeta}{\tan (w/2)}e^{-\frac{\zeta}{\tan^2 (w/2)}}-\tan^3(w/2) \left(1-e^{-\frac{\zeta}{\tan^2 (w/2)}}- \frac{\zeta}{\tan^2 (w/2)}\right) +
\\
&
+ 2\zeta\left[ \sinh \Big(\frac{\zeta}{2\tan^2 (w/2)}\Big)\cdot \left(\frac {\tan^2 (w/2)}{\zeta}-
\frac{z}{\tan (w/2)} \right)\right.+
\\
&\qquad\qquad\qquad\qquad\qquad\qquad
 \left. -\cosh\Big( \frac {\zeta}{2\tan^2 (w/2)}\Big)\right]\cdot \sum\limits_{j\geq k+1}p_j e^{-q_j+\eta/2} 
\end{array}
$$
$$
\begin{array}{rl}
\dot p_i = & p_i e^{-q_i+\eta/2}\left[ z \cosh\Big( \frac {\zeta}{2\tan^2 (w/2)}\Big)+ \tan(w/2) \sinh\Big(\frac {\zeta}{2\tan^2 (w/2)}\Big)\right] + 
\\
&+\sum\limits_{j\geq k+1} p_i p_j \sign(q_i-q_j) e^{-|q_i-q_j|}
\\
\dot q_i= & e^{-q_i+\eta/2}\left[ z \cosh\Big( \frac {\zeta}{2\tan^2 (w/2)}\Big)+ \tan(w/2) \sinh\Big(\frac {\zeta}{2\tan^2 (w/2)}\Big)\right]+
\\&+\sum\limits_{j\geq k+1} p_j e^{-|q_i-q_j|}
\end{array}
$$
which is a system of ODE with locally Lipschitz continuous right hand side that can be extended smoothly also at the value $w=\pi$. The initial data becomes
$$
z(\tau)=\lim_{t\to \tau^-} \sum_{i=1}^k p_i(t)\qquad w(\tau)=\pi \qquad  \eta(\tau)=2\bar q\qquad \zeta(\tau)= e_\tau 
$$
$$
p_i(\tau)=\lim_{t\to \tau^-}p_i(t)\qquad q_i(\tau)=\lim_{t\to \tau^-} q_i(t) \qquad i=k+1,\dots, N
$$
Thus there exists a unique solution of such a system  which provides a multi-peakon solution defined on some interval $[\tau,\tau'[$, up to the next interaction time.
As in Corollary of \cite[Section 7]{BF2}, once we prove that Camassa Holm equation is time reversible and the uniqueness of solutions of a Cauchy problem, we have that maximal number of peakon interaction is actually $k=2$, one with positive strength the other with negative one.

\section{A priori bounds}
This section is devoted to the study of some useful properties of the functions $u\in X_\alpha$. We start recalling an estimate for the $L^\infty-$norm of the $\H$ functions. We have
\begin{equation}
\label{sincichinequality}
\|f^2\|_{L^\infty}\leq \|f\|^2_{\H}
\end{equation}
This estimate give us a bound on the $L^\infty-$norm of the conservative solution $u$ of (\ref{prblCH}), in fact the conservation of the energy yields 
\begin{equation}
\label{normainfty}
\|u(t)\|_{L^\infty}\leq \|u(t)\|_\H=\sqrt{E^{\bar u}} \qquad \mbox{for every $t\geq 0$.}
\end{equation}691:
 Another fact is the behaviour of the functions $u\in X_\alpha$ at the infinitive. If we indicate with $C^{\alpha,u}$ the constant $\int_\R  (u^2+ u_x^2)e^{\alpha|x|} \,dx$, it holds
\begin{equation}
\label{uinfty}
\sup\limits_{x\in \R}\,u^2(x)e^{ \alpha |x|}\leq 2 C^{\alpha,u} \,.
\end{equation}
Indeed, the function
$$
f(x)\doteq u(t,x)e^{\frac\alpha 2|x| }. 
$$
belongs to $H^1(\R)$, moreover
$$
f_x=u_x e^{\frac\alpha 2|x| }+\frac \alpha 2 \sign(x)ue^{\frac\alpha 2|x| }
$$
and then, by using (\ref{sincichinequality}), we have 
$$
|f(x)|^2\leq\|f\|_\H^2\leq \int_\R [2u_x^2+ (1+\alpha)u^2]e^{\alpha |y|}\,dy \leq 2 C^{\alpha,u}\,.
$$
Now we study the behaviour at infinitive of the multi-peakon solutions of the Camassa Holm equation. 
\begin{lemma}({\rm A-priori} bounds)
Let $u$ be a multi-peakon solution to (\ref{prblCH}), with initial data $\bar u$ that satisfies (\ref{decayHP}). Then for every $t\in \R$ there exist a continuous function $C(t)$, which depends on $C^{\alpha, \bar u}$ and on the energy $E^{\bar u}$, such that
\begin{eqnarray}
\displaystyle
&&\int_{\R}[u^2(t,x)+u_x^2(t,x)]e^{\alpha|x|}\,dx\leq C(t)\,,
\label{weightedenergy}
\\&&
\displaystyle
\label{Pinfty}
\sup\limits_{x\in \R}\, \big|P_x^u(t,x)\big| e^{ \alpha |x|}\leq C(t)\,. 
\\&&
\|u_x\|_{L^1(\R)}\leq C(t)
\end{eqnarray} 
\end{lemma}
\begin{proof}
Since $|P_x^u|=P^u$, it is sufficient to prove the second inequality with $P^u$.
Setting 
$$
I(t)\doteq\int_{\R}[u^2(t,x)+u_x^2(t,x)]e^{\alpha|x|}\,dx\,,
$$
we want to achieve a differential inequality of the kind 
$$
\frac d{dt}I(t)\leq A+B\cdot I(t)\,,
$$
for some constants $A$ and $B$ which depends on the initial data $\bar u$. 
We start the discussion proving a preliminary estimate for the function $P^u$. By definition
\begin{equation}
\label{stimaP}
\int_{\R}P^u(t,x)e^{\alpha|x|}\,dx= \frac 12 \int_\R e^{\alpha |x|}\,dx 
\int_\R e^{-|x-y|}\left[u^2(t,y)+\frac{u_x^2(t,y)}{2}\right]dy 
\end{equation}
from this identity we can use Fubini's theorem to switch the order of the two integrals. Hence we compute the following integral 
\begin{equation}
\int_\R e^{\alpha|x|}e^{-|x-y|}\,dx=\frac {2\alpha}{1-\alpha^2}e^{-|y|}+\frac 2{1-\alpha^2}e^{\alpha |y|}\qquad \mbox{for every $y\in\R$}\,.
\label{expint}
\end{equation}
For future use, we observe that the equality (\ref{expint}) holds for $\alpha \in (-1,1)$. 
Substituting (\ref{expint}) in (\ref{stimaP}) and using the definition of the energy $E^{\bar u}$ we have
$$
\int_{\R}P^u(t,x)e^{\alpha|x|}dx\leq \frac{E^{\bar u}}{(1-\alpha^2)} +
\frac {1}{1-\alpha^2}\,I(t)\,.
$$
Having in mind the previous inequality, we are able to estimate the time derivative of the function $I$. From the equations (\ref{prblCH}) and (\ref{derevolution}) we have
$$
\begin{array}{rl}
\displaystyle
\frac d{dt}I(t)&
\displaystyle
\!\!=\!\! \int_{\R}\left[2u u_t+ \left(u_x^2\right)_t\right]\,e^{\alpha |x|}\,dx =
\int_{\R}\left[ -2u(u_x+P^u_x)+\frac 23 (u^3)_x- (uu_x^2)_x- 2 u_x P^u\right]\,e^{\alpha |x|}\,dx  
\\
&\displaystyle
\!\!\leq
-2\int_\R
(u P^u +u u_x^2)_x\,e^{\alpha|x|}\,dx \leq
\left. -2u (P^u + u_x^2)\,e^{\alpha|x|}\right|_{-\infty}^{\infty}
+ 2\alpha \int_\R |u| (P^u+ u_x^2)\,e^{\alpha|x|}\,dx
\\
&\displaystyle\!\!\leq to
\frac{2\alpha\|u\|_{L^\infty}}{1-\alpha^2}\left(E^{\bar u}+ 2 I\right)
\leq 
\frac{2\sqrt{E^{\bar u}}}{1-\alpha^2}\left[E^{\bar u}+ 2 I(t)\right]
\end{array}
$$
the previous inequality gives then a bound on the function $I$, that is
$$
I(t)\leq (C^{\alpha,\bar u}+E^{\bar u}/2)\exp\Big(\frac{4\sqrt{E^{\bar u}}}{1-\alpha^2}\, t \Big)\,.
$$691:
To achieve the estimate (\ref{Pinfty}), set
$$
K(t)\doteq \left\|P^u(t,\cdot)e^{\alpha|\cdot|}\right\|_{L^\infty}\,.
$$
Proceeding as before, fixed $x\in \R$ we compute the derivative w.r.t the time $t$ of the function $e^{-|x|}*{u^2_x}$. 
$$
\begin{array}{rl}
\displaystyle\frac \partial{\partial t} \left(e^{-|x|}*\frac{u^2_x}4\right)
&\displaystyle =\frac 14 \frac \partial{\partial t}\int_\R e^{-|x-y|}u_x^2(t,y)\,dy
\\
&\displaystyle
=\frac 12\int_\R e^{-|x-y|}\left[\Big( \frac {u^3}3- \frac{u u_x^2}2 - u P^u\Big)_x+ u P_x^u\right] dy
\\
&\displaystyle\leq \|u\|_{L^\infty} P^u+\frac  {\|u\|_{L^\infty}}2 \int_\R e^{-|x-y|}P^u(t,y)\,dy
\\
&\displaystyle\leq
 \|u\|_{L^\infty}P^u(t,x)+\frac{\|u\|_{L^\infty}}2 K(t)\int_\R e^{-\alpha|x|}e^{-|x-y|}\,dy
\\
&\leq \|u\|_{L^\infty} P^u(t,x)+\frac 1{1-\alpha^2}e^{-\alpha|x|}\|u\|_{L^{\infty}}K(t)
\end{array}
$$
in the same way, the derivative of $e^{-|x|}*u^2$ is
$$
\begin{array}{rl}
\frac \partial{\partial t}\left(e^{-|x|}*\frac{u^2}2\right)
&\displaystyle \leq \int_\R e^{-|x-y|}u \left(|P^u_x|+|u u_x|\right)\,dy
\\
&\displaystyle \leq\|u\|_{L^\infty} \left( 2 P^u(t,x)+\int_\R e^{-|x-y|}P^u(t,y)\,dy\right) 
\\
&\displaystyle \leq \|u\|_{L^\infty}\left( 2P^u(t,x)+\frac 1{1-\alpha^2} e^{-\alpha|x|}K(t)\right)\,.
\end{array}
$$
Multiplying the previous 691:two inequalities with $e^{\alpha|x|}$ we get
$$
\frac d{dt}K(t) \leq \left(3+\frac{ 2}{1-\alpha^2} \right) \sqrt{E^{\bar u}} K(t)
$$
which yields (\ref{Pinfty}). 

To achieve the last inequality, we write
\begin{eqnarray*}
\int_\R|u_x(y)|\,dy
&=&\!\!\!\!\int\limits_{ \{y :|u_x(y)|e^{\alpha|y|}<1 \} } \!\!\!\!|u_x(y)|\,dy
+\!\!\!\!\int \limits_{\{y:|u_x(y)|e^{\alpha|y|}>1\}} \!\!\!\!|u_x(y)|\,dy
\\
&\leq& \int_\R e^{-\alpha|y|}\,dy+\int_\R u_x^2(y)e^{\alpha|y|}\,dy\leq \frac 2\alpha + I(t)
\end{eqnarray*}
where the last estimate is given by (\ref{weightedenergy}).
\end{proof}
\section{Approximation of the initial data}
In this section we shall construct an approximation with an initial data with a multi-peakon function. Our aim is to approximate it with a sequence $u_\varepsilon$ which has an exponential decay at infinitive uniformly w.r.t $\varepsilon$.  
\begin{lemma}
Let $f\in X_\alpha$. Then for every $\varepsilon>0$ there exists a multi-peakon function $g$ of the form
$$
g(x)=\sum_{i=1}^N p_i e^{-|x-q_i|}
$$
such that 
\begin{eqnarray}
&&\|f-g\|_{\H}<\varepsilon
\\
&&\int_\R [g^2(x)+{g}_x^2(x)]e^{\alpha|x|}\,dx\leq C_0 
\end{eqnarray}
for some constant $C_0>0$691: which does not depend on $\eps$.
\end{lemma}
\begin{proof}
Let $\rho(x)\in \mathcal C^\infty_0$ be a cut-off function such that
\begin{itemize}
\item $\rho(x)\geq 0$
\item $\rho(x)=1$ for every $|x|\leq 1$, $\rho(x)=0$ for every $|x|> 2$
\item $\int_\R \rho(x)\,dx =1$ 
\end{itemize}
and $\rho_\eps(x)\doteq \frac 1\eps \rho(\frac x\eps) $ be a mollifiers sequence. Observe that for every $\eps>0$, $\tilde f(x)\doteq \rho_\eps*f(x)$ is a smooth function which approximates the function $f$ in $H^1-$norm
\begin{equation}
\|f-\tilde f\|_{\H}<C\eps
\label{apprcinf}
\end{equation}
 moreover it belongs to $X_\alpha$, indeed
\begin{eqnarray*}
\int_\R [\tilde f^2(x)+{\tilde f}_x^2(x)]e^{\alpha|x|}\,dx
&\leq& \int_\R \left[\int691:_\R (f^2(y)+f_x^2(y)) \rho_\eps(x-y) \,dy\right]e^{\alpha|x|}\,dx
\\
&\leq& \int_\R [f^2(y)+f_x^2(y)]\int_\R \rho_\eps(x-y) e^{\alpha|x|}\,dx \, dy
\\
&\leq& \int_\R [f^2(y)+f_x^2(y)] C_0 e^{\alpha|y|}\,dy =C_0 C^{\alpha,f}<\infty
\end{eqnarray*}
and $C_0$ is a constant which does not depend on $\eps$. From the previous inequality we can assert that for every $R>0$ one has
$\|\tilde f\|_{H^1(\R\setminus [-R,R])}\leq C_0 C_\alpha e^{-\alpha R}$ uniformly in $\eps>0$. We can choose thus $R_\eps$ big enough in order to have  
\begin{equation}
\|\tilde f\|_{H^1(\R\setminus [-R_\eps,R_\eps])}<\eps/2.
\label{fuori} 
\end{equation}

In the space $H^1([-R_\eps,R_\eps])$ we can approximate $f_\eps$ with a multi-peakon function.  By using the identity
$$
\frac 12\left(I-\frac{\partial^2}{{\partial x}^2}\right)e^{-|x|}=\delta_0
$$
the function $\tilde f$ can be rewritten in convolution form
$$ 
\tilde f= e^{-|x|}*\left(\frac{\tilde f - \tilde f_{xx}}2 \right)=\int_\R e^{-|x-y|}\cdot \frac{\tilde f(y) - \tilde f_{xx}(y)}2 \,dy\,.
$$
In the interval $[-R_\eps, R_\eps]$ the previous integral can now be approximated with a Riemann sum 
$$
g(x)= \sum_{i=-N}^N p_i e^{-|x-q_i|},\qquad 
\left\{
\begin{array}l
\dis q_i= \frac i N R_\eps 
\\
\dis p_i=\int_{q_{i-1}}^{q_i}\frac{\tilde f(y)-\tilde f_{xx}(y)}{2}\,dy\,.
\end{array}
\right.
$$
Choosing $N$ sufficiently large we obtain $\|\tilde f- g\|_{H^1([-R_\eps,R_\eps])}<\eps$. Together with (\ref{apprcinf}) and (\ref{fuori}) this last estimate yields the result.
\end{proof}

\section{Definition of the distance}\label{metricsection}
In this section we define a metric in order to control the distance between two solution of the equation (\ref{prblCH}). It is constructed as a Kantorovich-Wasserstein distance. Let $\T=[0,2\pi]$ the unit circle with the end points $0$ and $\pi$ identified. Consider the metric space $(\R^2\times \T,d^\diamondsuit)$, with distance 
$$
d^\diamondsuit((x,u,\omega),(x',u',\omega'))\doteq \min\{|x-x'|+|u-x'|+|\omega-\omega'|_*,1\}
$$
and for every function $u\in X_\alpha$, let define the Radon measure on $\R^2\times \T$
$$
\sigma^u(A)\doteq \int_{\{x\in \R: (x,u(x),2\arctan u_x(x))\in A\}} [1+u_x^2(x)]\,dx \qquad \mbox{for every Borel set $A$ of  $\R^2\times \T$}
$$
 The set $\mathcal F$ of transportation plans consists of the functions $\psi$ with the following properties:
\begin{enumerate}
\item \label{conduno} $\psi$ is absolutely continuous, is increasing and has an absolutely continuous inverse;
\item \label{condue} $\sup\limits_{x\in\R}|x-\psi(x)|e^{\alpha/2|x|}<\infty$;
\item \label{contre} $\int_{\R} |1-\psi'(x)|\,dx<\infty$.
\end{enumerate}
\psfrag{x}{$x$}
\psfrag{u}{$u$}
\psfrag{v}{$v$}
\begin{figure}[ht]
\centerline{
\includegraphics[width=10cm]{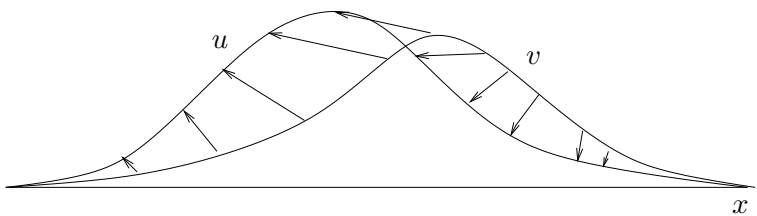}
}\caption{Transportation plan.\label{ftmass}}
\end{figure}  
The conditions \ref{condue} and \ref{contre} are not restrictive. Indeed, thanks to the exponential decay of functions $u,v \in X_\alpha$, the measures $\sigma^u$ and $\sigma^v$ located on the graph of $u$ and $v$ respectively, have small mass at the infinitive, and then a transportation plan which transports mass from one to the other can be almost the identity $\psi(x) \approx x$ (see fig. \ref{ftmass}).
In order to define a distance in the space $X_\alpha$, we consider an optimization problem over all possible transportation plans. Given two functions $u,\,v$ in $X_\alpha$, we introduce two further measurable functions, related to a transportation plan $\psi$:
\begin{eqnarray}
&&
\label{non-phi1}\phi_1(x)\doteq \sup \big\{ \theta \in [0,1] \st \theta \cdot(1+u_x^2(x))\leq \left(1+v_x^2(\psi(x))\right)\psi'(x)\big\},
\\
&&\label{non-phi2}\phi_2(x)\doteq \sup \big\{ \theta \in [0,1] \st  1+u_x^2(x)\leq \theta\cdot \left(1+v_x^2(\psi(x))\right)\psi'(x)\big\}.
\end{eqnarray}
The functions $\phi_1, \phi_2$ can be seen as weights that take into account the difference of the masses of the measure $\sigma^u$ and $\sigma^v$. In fact, from the definitions (\ref{non-phi1})-(\ref{non-phi2}) one has
$$
\phi_1(x) (1+u_x^2(x))=\phi_2(\psi(x)) (1+v_x^2(\psi(x)))\psi'(x)\qquad \mbox{for a.e. $x\in\R$}.
$$
According to the definitions, the identity $\max\{\phi_1(x),\phi_2(x)\}\equiv 1$ holds.
Altough the two measures $\phi_1\sigma^u$ and $\phi_2\sigma^v$ have not finite mass, they satisfy $\phi_1\sigma^u(A)=\phi_2\sigma^v(A)$ for every bounded Borel set $A\subset \R^2\times \mathbb T$. Thus, the functions $\phi_1$ and $\phi_2$ represent the percentage of mass actually transported from one measure to the other. 
A distance between the two functions $u,\,v$ in $X_\alpha$ can be characterized in the following way.

For every $\psi\in \mathcal F$, let ${\bf X}^u=(x,u(x),2\arctan u_x(x))$ and ${\bf X}^v=(\psi(x),v(\psi(x)),2\arctan v_x(\psi(x)))$ and consider the functional
$$
J^\psi(u,v)=\int_\R d^\diamondsuit ({\bf X}^u,{\bf X}^v)\phi_1(x)(1+u_x^2(x))\,dx+ \int_\R \left|1+u_x^2(x)-(1+v_x^2(\psi(x)))\psi'(x)	\right|\,dx\,.
$$
Since the previous function is well defined for every $\psi\in \mathcal F$, we can define 
$$
J(u,v)\doteq \inf_{\psi \in \mathcal F} J^\psi (u,v).
$$
The function $J$ here defined is thus a metric on the space $X_\alpha$ (see \cite{BF2}). 

\medskip
\section{Comparison with other topologies}
\begin{lemma}\label{L1J}
For every $u,\,v\in X_\alpha$ one has 
\begin{equation}
\label{h1el1}
\frac 1C\cdot \|u-v\|_{L^1(\R)}\leq J(u,v)\leq C\cdot \|u-v\|_{\H}.
\end{equation}
Let $(u_n)$ be a Cauchy sequence for the distance $J$ such that $C^{\alpha, u_n}\leq C_0$ for every $n\in \N$. Then
\begin{enumerate}
\item[i] There exists a limit function $u\in X_\alpha$ such that $u_n\to u$ in $L^\infty$ and the sequence of derivatives ${u_n}_x$ converges to $u_x$ in $L^p(\R)$ for $p \in [1,2[$.
\item[ii] Let $\mu_n$ the absolutely continuous measure having density ${u_n}_x^2$ w.r.t. Lebesgue measure. then one has the weak convergence $\mu_n\rightharpoonup \mu$ for some measure $\mu$ whose absolutely continuous part has density $u_x^2$. 
\end{enumerate}
\end{lemma}
\begin{proof}
The first inequality of (\ref{h1el1}) can be achieved by estimating the area between the two functions $u$ and $v$. For every $\psi\in \F$ we can write
$$
\int_\R |u-v|\, dx =\left( \int_{S_1}+\int_{S_2}\right) |u-v|\,dx
$$ 
where the two subsets $S_1$ and $S_2$ are
\begin{itemize}
\item $S_1= \{x: |x-\psi(x)|\leq~1~\}=\cup_j [x_{2j-1},x_{2j}]$, where in this union we have to take into account that these intervals may be either finite or infinite, possibly having $x_j=\pm \infty$ for some $j$, 
\item $S_2=\{x: |x-\psi(x)|>1\}$. 
\end{itemize}
The integral over $S_2$ can be estimate in the following way:
\begin{equation}
\label{psigrosso}
\int_{S_2}|u(x)-v(x)|\,dx\leq (\|u\|_{L^\infty}+\|v\|_{L^\infty}) \int_\R|x-\psi(x)|\,dx\leq (E^{\bar u}+ E^{\bar v})J(u,v).
\end{equation}
The last estimate is given by the definition of the functional $J$.
\psfrag{An}{$A_{n}$}
\psfrag{An+1}{$A_{n+1}$}
\psfrag{An+k}{$A_{n+k-1}$}
\psfrag{An+j}{$A_{n+k}$}
\psfrag{u}{$u$}
\psfrag{v}{$v$}
\psfrag{x0}{$x_0$}
\psfrag{x1}{$x_1$}
\psfrag{x2k}{$x_{2k}$}
\psfrag{x2k+1}{$x_{2k+1}$}
\psfrag{xk}{$x_k$}
\psfrag{x-psix<1}{$S_1$}
\psfrag{x-psix>1}{$S_2$}
\psfrag{n}{$n$}
\psfrag{n+1}{$n+1$}
\psfrag{n+j}{$n+j$}
\psfrag{n+k}{$n+k$}
\psfrag{...}{$...$}
\psfrag{PQ}{$\overline {Q_u(n)\,Q_v(n)}$}
\begin{figure}[ht]
\centerline{\includegraphics[width=10cm]{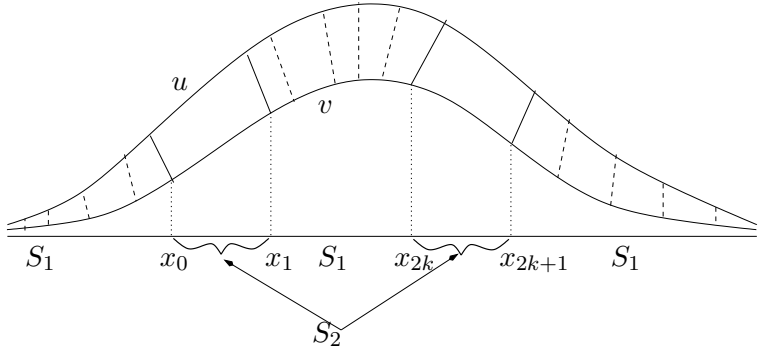}}
\caption{$L^1-$distance between two functions.\label{L1norma}}
\end{figure}

As far as the integral over $S_1$ is concerned, the integral over $S_1$ can be viewed as a sum of the area of the regions $A_j$ in the plane $\R^2$, bounded by the graph of the curves $u$, $v$ and by the segments with slope $\pm 1$ that join the points $Q_u(x_{2j-1})=(x_{2j-1},u(x_{2j-1}))$ and $Q_v(x_{2j})=(\psi(x_{2j}),v(\psi(x_{2j})))$, where $\{x_i\}=\partial S_1$. 
We have
$$
\int_{S_1}|u(x)-v(x)|\,dx\leq \sum_j \meas(A_j)\,.
$$
The measure of the subset $A_j$ is the area sweeped by the segment $\overline {Q_u(x)\,Q_v(x)}$. Recalling that in every set $A_j$ the function $\psi$ satisfies $ |x|-1\leq |\psi(x)|\leq |x|+1$, a bound on this area is given by
$$
\meas (A_j) \leq
\displaystyle
\int_{x_{2j-1}}^{x_{2j}} (|x-\psi(x)|+|u(x)-v(\psi)|)[1+u_x^2+(1+v_x^2(\psi ))\psi']\,dx
$$
and then 
$$
\begin{array}{rl}
\dis
\int_{S_1}|u(x)-v(x)|\,dx
&\dis 
\leq \int_{S_1}(|x-\psi(x)|+|u(x)-v(\psi)|)[1+u_x^2+(1+v_x^2(\psi ))\psi']\,dx
\\
&\dis\leq J^\psi(u,v) + J^{\psi^{-1}} (u,v)
\end{array}
$$
this inequality, together with (\ref{psigrosso}), yields to 
$$
A\leq C(\bar u,\bar v)  J(u,v).
$$
The proof of the second part of the lemma is perfectly similar to the one of the periodic case, once we take into account the exponential decay of the sequence $u_n$.   
\end{proof}
\section{Stability of solutions w.r.t initial data}\label{stabcor}
Let $u_0$ and $v_0$ be two multi-peakon initial data. The technique developed in Section \ref{ODEsystem} ensures the existence of two multi-peakon solution $u(t),v(t)$ for (\ref{decayHP}) which conserve the energy unless interaction of peakon occurs. Suppose then that within a given interval $[0,T]$ no interaction occurs neither for $u(t)$ nor for $v(t)$. 
The aim of this section is to prove the continuity of the functional $J$ w.r.t. the initial data, more clearly we prove that there exists a continuous, positive function $C(t)$ such that for $t\in [0,T]$ one has
$$
J(u(t),v(t))\leq C(t) J(u_0,v_0).
$$
\begin{lemma}
If $u(t)$ and $v(t)$ are two multi-peakon solutions defined in the interval $[0,T]$ in which no interaction occurs, then there exists a positive, continuous function $c(t)$ which depends only on the energies $E^u$, $E^v$ of the two solutions, such that
\begin{equation}
\frac d{dt} J(u(t),v(t))\leq c(t) J(u(t),v(t)) \qquad \mbox{for all $t\in[0,T]$}. 
\end{equation}
\end{lemma}
\begin{proof}
We compute the time derivative of the function $J^\psi(u(t),v(t))$ with a particular choice of the transportation plan $\psi=\psi_{(t)}$. Given any $\psi_0\in \mathcal F$, at every time $t\in [0,T]$ we construct $\psi_{(t)}$ by transporting the function $\psi_0$ along the characteristic curves. More precisely, since no interaction between peakon occurs in the interval $[0,T]$, the functions $u(t,\cdot),\,v(t,\cdot)$ are Lipschitz continuous, then the flows $\vphi^t_u$, $\vphi^t_v$ solutions of the Cauchy problems
\begin{eqnarray*}
&&\frac d{dt} \vphi^t_u(x)=u(t,\vphi_u^t(x))\qquad \vphi^0_u(x)=x,
\\
&&\frac d{dt}\vphi^t_v(y)=v(t,\vphi_v^t(y))\qquad \vphi^0_v(y)=y,
\end{eqnarray*}
which are the characteristics curves associated to the equation (\ref{decayHP}), are well defined. Now, let $x\in \R$. $\psi_{(t)}$ is defined as the composition
\begin{equation}
\psi_{(t)}(x)\doteq \vphi_v^t\circ \psi_0\circ \left(\vphi^t_u \right)^{-1}(x),
\end{equation}
that is 
$$
\psi_{(t)}(\vphi_u^t(y))=\vphi_v^t(\psi_0(y)).
$$
The function $\psi_{(t)}$ belongs to $\mathcal F$, and hence $J^{\psi_{(t)}}$ is well defined, in fact
\begin{enumerate}
\item By the property \ref{conduno} of the function $\psi_0$ and uniqueness of solution of ODE, the function $\psi_{(t)}$ is an  increasing function.
\item Let $x\in \R$ and $\vphi_u^t(y)$ be the characteristic curve passing through $x$ at time $t$. Evaluating $|x-\psi_ {(t)}(x)|e^{\alpha/2|x|}$ along this characteristic curve, and computing the derivative w.r.t. $t$ we obtain
$$
\frac d{dt}|\vphi_u^t(y)-\vphi_v^t(\psi_0(y))|e^{\alpha/2|\vphi_u^t(y)|}\leq 
\left[|u(t,x)-v(t,\psi_{(t)}(x))|+\frac \alpha 2 |u(t,x)|\cdot |x-\psi_{(t)}(x)|\right]e^{\alpha/2|x|}
$$ 
by properties (\ref{uinfty}), (\ref{weightedenergy}), and since $u,\,v$ are Lipschitz continuous in $[0,T]$, there exists two  $L^\infty$ functions $c_1(t)$, $c_2(t)$ such that 
$$
\frac d{dt}|x-\psi_{(t)}(x)|e^{\alpha/2|x|}\leq c_1(t)|x-\psi_{(t)}(x)|e^{\alpha/2|x|}+c_2(t)
$$ 
by Gronwall Lemma and the hypothesis $|x-\psi_0(x)|e^{\alpha/2|x|}\leq C_0$, the previous inequality gives the property \ref{condue} for $\psi_{(t)}$
\begin{equation}
\label{stpsi}
|x-\psi_{(t)}(x)|e^{\alpha/2|x|}\leq C_1(t)\doteq \left(C_0+ \int_0^t c_2(s)\,ds\right)e^{\int_0^t c_1(s)\,ds}
\end{equation}
\item The last property can be proved by changing the integration variable $x= \vphi_u^t(y)$ 
$$
\begin{array}{rl}
\dis 
\int_\R |1-\psi_{(t)}(x)|\,dx
&
\dis 
=\int_\R |1-\psi_{(t)}(x)|(\vphi_u^t)'(y)\,dy = \int_\R |(\vphi_u^t)'(y)-(\vphi_v^t)'(\psi_0(y))\psi_0'(y)|\,dy
\\
&\leq \dis \int_\R |(\vphi_u^t)'(y)-1|\,dy +\int_\R |(\vphi_v^t)'(y)-1|\,dy+\int_\R|1-\psi_0'(y)|\,dy.
\end{array}
$$
Since
$$
|(\vphi_u^t)'(y)-1|\leq \int_0^t |u_x(s,x)|\cdot |(\vphi_u^s)'(y)-1|\,ds+\int_0^t |u_x(s,x)|\,ds
$$
(and a similar estimate for $\vphi_v^t$) and $u_x,v_x\in L^\infty$, by the Gronwall lemma the first two integrals of the previous formula are bounded by an absolutely continuous function $C(t)$ in the interval $[0,T]$ and then also property \ref{contre} holds.
\end{enumerate}
At the transportation plan $\psi_{(t)}$ we associate the functions $\phi_1^{(t)},\,\phi_2^{(t)}$ defined according to  (\ref{non-phi1}), (\ref{non-phi2}), the functional $J^{\psi_{(t)}}$ is thus 
$$
J^{\psi_{(t)}}(u(t),v(t))=\!
\int_\R d^\diamondsuit(\mathbf X^u(t),\mathbf X^v(t))\phi_1^{(t)}(x)(1+u_x^2(x))dx+
\int_\R \left|1+u_x^2(x)
-(1+v_x^2(\psi_{(t)}(x)))\psi_{(t)}'(x)\right|dx.
$$
By deriving $J^{\psi_{(t)}}(u(t),v(t))$ w.r.t. $t$ and computing the change of variables along the characteristics, the previous derivative can be estimate by the sum of the following terms (we leave out the dependence on the integrable variable when it is not essential)
\begin{itemize}
\item $\dis
I_1=\int_\R |u(t,x)-v(t,\psi_{(t)}(x))|\phi_1^{(t)}(x)(1+u_x^2(t,x))\,dx\leq (1+\|u(t)\|_{L^\infty}+\|v(t)\|_{L^\infty}) J^{\psi_{(t)}}(u(t),v(t))\,,
$
\item $
\dis I_2=\int_\R|P^u_x(t,x)-P^v_x(t,\psi_{(t)}(x))|\phi_1^{(t)}(x)(1+u_x^2(t,x))\,dx\,,
$ 
\item $I_3=
\dis
\int_\R \left|\frac {2u^2(t)-u_x^2(t)-2P^u(t)}{1+u_x^2(t)}\right.
\left.-\frac{2v^2(t,\psi_{(t)})-v_x^2(t,\psi_{(t)})-2P^v(t,\psi_{(t)})} {1+v_x^2(t,\psi_{(t)})}\right|\cdot
\phi_1^{(t)} (1+u_x^2(t))\,dx\,,
$
\item the term due to the variation of the base measure 
$$
\begin{array}{rl}
\dis
I_4=2\int_\R d^\diamondsuit(\mathbf X^u(t),\mathbf X^v(t))\cdot u_x(t)(u^2(t)- P^u(t))\,dx\,,
\end{array}
$$
\item and the terms due to the variation of the excess mass
$$
I_5=\frac d{dt}\int_\R\left|1+u_x^2(t)-(1+v_x^2(t,\psi_{(t)}))\psi_{(t)}'\right|\,dx\,.
$$
\end{itemize}
Let us start to estimate the term $I_2$. By definition, the difference of $P^u$ and $P^v$ is written in convolution form
$$
\begin{array}{rl}
 &\dis\left| \int_\R\left \{  e^{-|x-y|}\sign(x-y)\left[u^2(t,y)+\frac{u^2_x(t,y)}{2}\right]\,dy \right.\right.
\\
&\left.\left.\dis \qquad\qquad -e^{-|\psi_{(t)}(x)-\psi_{(t)}(y)|}\sign(\psi_{(t)}(x)-\psi_{(t)}(y))\left[v^2(t,\psi_{(t)}(y))+\frac{v^2_x(t,\psi_{(t)}(y))}{2}\right]\psi_{(t)}'(y)\right\}\,dy\right|
\end{array}
$$
then in $I_2$ appear the following integrals
\begin{eqnarray*}
\displaystyle
&&A=\int_\R(1+u_x^2(t,x))\int_\R e^{-|x-y|}|u^2(t,y)-v^2(t,y)|\,dy\,dx
\\
&&B=\int_\R(1+u_x^2(t,x))\left|\int_\R e^{-|x-y|}\sign(x-y)[v^2(t,y)-v^2(t,\psi_{(t)}(y))\psi_{(t)}'(y)]\,dy \,\right|\,dx
\\
&&C=\int_\R(1+u_x^2(t,x))\int_\R\left|e^{-|x-y|}\sign(x-y)-e^{-|\psi_{(t)}(x)-\psi_{(t)}(y)|}\sign(\psi_{(t)}(x)-\psi_{(t)}(y))\right| \cdot
\\
&&\qquad\qquad\qquad\qquad\qquad\qquad
\cdot\left[v^2(t,\psi_{(t)} (y))+\frac {v_x^2(t,\psi_{(t)} (y))}2\right]\psi_{(t)}'(y)\,dy\,dx
\\
&&D=\frac 12\int_\R(1+u_x^2(t,x)) \left|\int_\R e^{-|x-y|}\sign(x-y)[u_x^2(t,y)-v_x^2(t,\psi_{(t)}(y))\psi_{(t)}'(y)]\,dy\,\right|\,dx
\end{eqnarray*}

\bigskip\noindent{\bf A.} Switching the order of the two integrals, the term $A$ is bounded by the $L^1$-norm of the difference between $u$ and $v$:
$$
\begin{array}{rl}
A
& 
\dis 
\leq(\|u\|_{L^\infty}+\|v\|_{L^\infty})\int_\R |u(t,y)-v(t,y)|\int_\R e^{-|x-y|}(1+u_x^2(t,x))\,dx\,dy
\\
&\dis
\leq (2+E^u)(\|u\|_{L^\infty}+\|v\|_{L^\infty})\|u(t)-v(t)\|_{L^1}
\end{array}
$$
and then, by Lemma \ref{L1J}, $A\leq C(\bar u,\bar v)J(u(t),v(t))$.

\bigskip
\noindent{\bf B.}
Define 
$$
F(y)\doteq\int_{-\infty}^y (v^2(z)-v^2(\psi_{(t)}(z))\psi_{(t)}'(z))\,dz=\int_{\psi_{(t)}(y)}^{y}v^2(z)\,dz
$$
we have,
integrating by parts
$$
\begin{array}{rl}
\dis
\left|\int_\R e^{-|x-y|}\sign(x-y)F'(y)\,dy\,\right|
&
\dis
\leq 2|F(x)|+\int_R e^{-|x-y|}|F(y)|\,dy
\\
&\dis
\leq \|v\|^2_{L^\infty}|x-\psi(x)|+\int_\R e^{-|x-y|}|y-\psi(y)|\,dy
\end{array}
$$
moreover, substituting the previous expression into the term $B$ we obtain
\begin{eqnarray*}
B&\!\!\!\!\!\!\!\!\!\!&\leq 2E^{\bar v}J^{\psi_{(t)}}(u(t),v(t))+\int_\R (1+u_x^2(t,x)) \int_\R e^{-|x-y|} |y-\psi_{(t)}(y)|\,dy\,dx
\\
&\!\!\!\!\!\!\!\!\!\!&=2E^{\bar v}\left\{J^{\psi_{(t)}}(u(t),v(t))+\int_\R |y-\psi_{(t)}(y)|\int_\R(1+u_x^2(t,x))e^{-|x-y|} \,dx\,dy\right\}
\\
&\!\!\!\!\!\!\!\!\!\!&\leq 2E^{\bar v}(3+E^{\bar u})\cdot J^{\psi{(t)}}(u(t),v(t)).
\end{eqnarray*}

\bigskip
\noindent {\bf C.} Observe that since the function $y\mapsto \psi_{(t)}(y)$ is non decreasing, the quantities $x-y$ and $\psi_{(t)}(x)-\psi_{(t)}(y)$ have the same sign, and since the function $t \mapsto e^{-|t|}$ is Lipschitz continuous either in $(-\infty,0)$ or in $(0,+\infty)$ we have
$$
\begin{array}{rl}
\dis
\left|e^{-|x-y|}-e^{-|\psi_{(t)}(x)-\psi_{(t)}(y)|}\right|
&\leq \dis
e^{-\min\{|x-y|,|\psi_{(t)}(x)-\psi_{(t)}(y)|\}}\left||x-y|-|\psi_{(t)}(x)-\psi_{(t)}(y)|\right|
\\
&
\dis
\leq
e^{-\min\{|x-y|,|\psi_{(t)}(x)-\psi_{(t)}(y)|\}}\left(|x-\psi_{(t)}(x)|+|y-\psi_{(t)}(y)|\right)
\end{array}
$$ 
now 
$$
-\min\{|x-y|,|\psi_{(t)}(x)-\psi_{(t)}(y)|\}\leq -|x-y|+2C_1(t),
$$
where $C_1(t)$ is the function (\ref{stpsi}), related to the property \ref{condue} of $\psi_{(t)}$, then
$$
\begin{array}{rl}
C
\leq 
&\!\!\dis e^{C_1(t)}\!\!\int_\R |y-\psi(y)|
\left[v^2(t,\psi_{(t)}(y)) + \frac {v_x^2(t,\psi_{(t)}(y))}2\right]\psi_{(t)}'(y)\cdot\int_\R(1+u_x^2(t,x))e^{-|x-y|}\,dx\,dy+
\\
&+\dis 2 E^{\bar v}\int_\R(1+u_x^2(t,x))|x-\psi_{(t)}(x)|\,dx
\\
\leq 
& \!\!\dis
\left[(2+E^{u})e^{C_1(t)}(1+\|u\|_{L^\infty})+ 2 E^{\bar v}\right] J^{\psi_{(t)}}(u(t),v(t))
\end{array}
$$

\bigskip
\noindent{\bf D.} 
Here we can use the estimate given by the change in base measure. Since  
$$
\int_\R \Big|1+u_x^2(t,x)- \big(1+v_x^2(t,\psi_{(t)}(x))\big)\psi_{(t)}'(x)\Big|\, dx\leq J^{\psi_{(t)}}(u(t),v(t))\,,
$$
we obtain
$$
\begin{array}{rl}
D \leq&
\dis
\frac 12 \int_\R(1+u_x^2(t,x)) \int_\R e^{-|x-y|}
\left|(1+u_x^2(t,y))-(1+v_x^2(t,\psi_{(t)}(y)))\psi_{(t)}'(y)\right|  \,dy\,\,dx
\\
\quad&
\dis  + \frac 12 \int_\R(1+u_x^2(t,x))\cdot 
\left|\int_\R  e^{-|x-y|}\sign(x-y) \left[\psi_{(t)}'(y)-1\right]dy\, \right|\,dx
\\
\leq& 
\dis
(1+ E^{\bar u})J^{\psi_{(t)}}(u(t),v(t))+
 \int_\R(1+u_x^2(t,x))\left|(\psi_{(t)}(x)-x)-\frac{e^{-|x-y|}}2\int_\R (\psi_{(t)}(y)-y)\,dy\right|\,dx 
\\
\leq &
2(2+ E^{\bar u})J^{\psi_{(t)}}(u(t),v(t))
\dis
\end{array}
$$
where in the last estimate we integrated by part as in the term $B$. 

\bigskip
\noindent The control for the terms $I_3$, $I_4$ and $I_5$ can be obtained exactly as the ones in \cite{BF2}, whom we refer the reader to.  
The previous estimates implies that there exists a smooth function $C=C^{\bar u,\bar v}(t)$ which depends only to the variable $t$ and to the initial data $\bar u,\ \bar v$ such that 
$$
\frac{d}{dt} J^\psi(u,v)\leq C^{\bar u,\bar v}(t) J^{\psi}(u,v)
$$
which yields
$$
J(u(t),v(t))\leq J(u(s), v(s)) e^{\left|\int_s^t C^{\bar u,\bar v}(\sigma)\,d\sigma\right|}\qquad {\mbox {for every $s,t\in\R$.}}
$$
\end{proof}

\end{document}